\RequirePackage{ifpdf}
\ifpdf 
\documentclass[pdftex]{sigma}
\else
\documentclass{sigma}
\fi

\begin{document}

\allowdisplaybreaks

\renewcommand{\thefootnote}{$\star$}

\renewcommand{\PaperNumber}{079}

\FirstPageHeading

\ShortArticleName{Non-Gatherable Triples for Non-Aaf\/f\/ine Root Systems}

\ArticleName{Non-Gatherable Triples for Non-Af\/f\/ine Root Systems\footnote{This paper is a
contribution to the Special Issue on Kac--Moody Algebras and Applications. The
full collection is available at
\href{http://www.emis.de/journals/SIGMA/Kac-Moody_algebras.html}{http://www.emis.de/journals/SIGMA/Kac-Moody{\_}algebras.html}}}

\Author{Ivan CHEREDNIK and Keith SCHNEIDER}

\AuthorNameForHeading{I.~Cherednik and K.~Schneider}

\Address{Department of Mathematics, UNC
Chapel Hill, North Carolina 27599, USA}
\Email{\href{mailto:chered@email.unc.edu}{chered@email.unc.edu}, \href{mailto:schneidk@email.unc.edu}{schneidk@email.unc.edu}}

\ArticleDates{Received September 03, 2008, in f\/inal form November 08,
2008; Published online November 14, 2008}

\Abstract{This paper contains a complete description of minimal non-gatherable
triangle triples in the lambda-sequences for the classical root systems,
$F_4$ and $E_6$. Such sequences are associated with reduced decompositions
(words) in af\/f\/ine and non-af\/f\/ine Weyl groups. The existence of the
non-gatherable triples is a combinatorial obstacle for using the
technique of intertwiners for an explicit description of the irreducible
representations of the (double) af\/f\/ine Hecke algebras, complementary
to their algebraic-geometric theory.}


\Classification{20H15; 20F55}

\newcommand{\comment}[1]{}
\renewcommand{\tilde}{\widetilde}
\renewcommand{\hat}{\widehat}

\renewcommand{\tilde}{\widetilde}
\renewcommand{\hat}{\widehat}

\newcommand{\BR}{{\mathbb R}}
\newcommand{\BQ}{{\mathbb Q}}
\newcommand{\BC}{{\mathbb C}}
\newcommand{\BP}{{\mathbb P}}
\newcommand{\BZ}{{\mathbb Z}}
\newcommand{\BN}{{\mathbb N}}
\newcommand{\BS}{{\mathbb S}}

\newcommand{\cH}{{\mathcal H}}
\newcommand{\cA}{{\mathcal A}}
\newcommand{\cB}{{\mathcal B}}
\newcommand{\ccF}{{\mathfrak F}}
\newcommand{\cD}{{\mathcal D}}
\newcommand{\cL}{{\mathcal L}}
\newcommand{\cF}{{\mathcal F}}
\newcommand{\cP}{{\mathcal P}}
\newcommand{\cX}{{\mathcal X}}
\newcommand{\cY}{{\mathcal Y}}
\newcommand{\cS}{{\mathcal S}}
\newcommand{\cSol}{\hbox{$\mathcal Sol$}}
\newcommand{\cT}{\hbox{$\mathcal T$}}

\newcommand{\Z}{{\mathbb Z}}
\newcommand{\Q}{{\mathbb Q}}
\newcommand{\N}{{\mathbb N}}
\newcommand{\C}{{\mathbb C}}
\newcommand{\R}{{\mathbb R}}

\newcommand{\CH}{{\mathcal H}}
\newcommand{\CA}{{\mathcal A}}

\def\HH{\mbox{${\mathcal H}$\kern-5.2pt${\mathcal H}$}}


\def\der{\partial}
\def\tensor{\otimes}
\def\gam{\gamma} \def\Gam{\Gamma}
\def\del{\delta} \def\Del{\Delta}
\def\kap{\kappa}
\def\lam{\lambda} \def\Lam{\Lambda}
\def\Comp{{\mathbb C}}
\def\sM{{\mathcal M}}

\comment{
\newtheorem{theorem}{Theorem}[section]
\newtheorem{proposition}[theorem]{Proposition}
\newtheorem{definition}[theorem]{Definition}
\newtheorem{lemma}[theorem]{Lemma}
\newtheorem{corollary}[theorem]{Corollary}
\newtheorem{notation}[theorem]{Notation}
\newtheorem{remark}[theorem]{Remark}
\newtheorem{example}[theorem]{Example}
\newtheorem{note}{Note}
}

\newtheorem{maintheorem}[theorem]{Main Theorem}
\newtheorem{theorem }{Theorem}[section]
\newtheorem{maintheorem }[theorem]{Main Theorem}
\newtheorem{proposition }[theorem]{Proposition}
\newtheorem{definition }[theorem]{Definition}
\newtheorem{lemma }[theorem]{Lemma}
\newtheorem{corollary }[theorem]{Corollary}
\newtheorem{notation }[theorem]{Notation}
\newtheorem{remark }[theorem]{Remark}
\newtheorem{example }[theorem]{Example}

\newtheorem{ maintheorem }[theorem]{Main Theorem}
\newtheorem{ theorem}{Theorem}[section]
\newtheorem{ proposition}[theorem]{Proposition}
\newtheorem{ definition}[theorem]{Definition}
\newtheorem{ lemma}[theorem]{Lemma}
\newtheorem{ corollary}[theorem]{Corollary}
\newtheorem{ notation}[theorem]{Notation}
\newtheorem{ remark}[theorem]{Remark}
\newtheorem{ example}[theorem]{Example}

\newtheorem{thm}{Theorem}[section]
\newtheorem{prop}[thm]{Proposition}
\newtheorem{lem}[thm]{Lemma}
\newtheorem{cor}[thm]{Corollary}
\newtheorem{conj}[thm]{Conjecture}
\newtheorem{con}[thm]{Conjecture}
\newtheorem{dfn}[thm]{Definition}
\newtheorem{df}[thm]{Definition}
 \newcommand{\rem}{{\bf Comment.\ }}
 \newcommand{\rmk}{{\bf Comment.\ }}
 \newcommand{\exmp}{{\bf Example.\ }}
 \newcommand{\ex}{{\bf Example.\ }}
 \newcommand{\prob}{{\bf Problem.\ }}

\renewcommand{\thenote}{}
\newtheorem*{acka}{Acknowledgments}
\newtheorem{ack}{Acknowledgments}
\renewcommand{\theack}{}
\renewcommand{\appendixname}{\bf Appendix}

\hyphenation{
ap-pen-dix as-ymp-tot-ic at-trib-uted at-trib-ut-able
Bry-li-n-sky com-mu-ta-tion de-ge-ne-rate
de-riv-a-tive dis-trib-ute equi-vari-ant ex-tra-or-di-nary
geo-met-ric griev-ance griev-ous grad-ed ho-lo-no-my ho-mo-thetic
in-fin-ite-ly in-fin-i-tes-i-mal Ha-rish Cha-n-dra mul-ti-plic-able
non-euclid-ean non-iso-mor-phic non-smooth par-a-digm
par-a-bol-ic pa-rab-o-loid pa-ram-e-trize phe-nom-e-non
post-script pseu-do-dif-fer-en-tial pseu-do-fi-nite
qua-drat-ics quad-ra-ture Han-kel rec-tan-gle semi-def-i-nite
set-up wide-spread Euler-ian Feb-ru-ary Gauss-ian Grothen-dieck
Hamil-ton-ian Her-mi-t-ian her-mi-t-ian Jan-u-ary
Japan-ese Ka-shi-wa-ra Kor-te-weg Le-gendre No-vem-ber Rie-mann-ian
Sep-tem-ber Za-mo-lo-d-chi-kov Kni-zh-nik quan-tum Op-dam
Mac-do-nald Ca-lo-ge-ro Su-ther-land Mo-ser
Ol-sha-net-sky  Pe-re-lo-mov in-de-pen-dent ope-ra-tors
cy-clo-to-mic ra-tio-nal de-gen-er-a-tion
in-ter-est-ing de-for-ma-tions de-for-ma-tion pro-ce-dure
fol-lows ope-ra-tors  pre-serve suf-fices ap-proach
for-mu-las con-sider its com-ple-tion cor-re-spond-ing
au-to-mor-phism be-cause pro-por-tional fi-nal-ly let-ting
equi-v-a-lence ge-n-er-al-ized Mac-do-nald iden-ti-ties
cor-re-s-pond sub-dia-grams par-ti-tion na-t-u-ral-ly
or-dered stan-dard de-for-ma-tion ar-gu-ment com-bined
sphe-r-i-cal rep-re-sen-ta-tions tri-go-no-me-t-ric
ge-n-er-al-ly speak-ing pri-m-it-ive ir-re-du-cible
sum-ma-tion  rep-re-sen-ta-tives pro-por-ti-o-na-li-ty
ultra-sphe-ri-cal Ro-gers}

\def\ffor{\quad\hbox{ for }\quad}
\def\wwhen{\quad\hbox{ when }\quad}
\def\wwhere{\quad\hbox{ where }\quad}
\def\aand{\quad\hbox{ and }\quad}
\def\for{\  \hbox{ for } \ }
\def\iif{ \ \hbox{ if } \ }
\def\when{ \ \hbox{ when } \ }
\def\where{\  \hbox{ where } \ }
\def\and{\  \hbox{ and } \ }
\def\and{\  \hbox{ and } \ }
\def\oor{\  \hbox{ or } \ }
\def\equal{\stackrel{\rm def}{= \kern-3pt =}}

\def\la{\lambda}
\def\La{\Lambda}
\def\om{\omega}
\def\Om{\Omega}
\def\Th{\Theta}
\def\th{\theta}
\def\al{\alpha}
\def\be{\beta}
\def\ga{\gamma}
\def\ep{\epsilon}
\def\up{\upsilon}
\def\Up{\Upsilon}
\def\de{\delta}
\def\De{\Delta}
\def\ka{\kappa}
\def\kapp{\hbox{\bf \ae}}
\def\si{\sigma}
\def\Si{\Sigma}
\def\Ga{\Gamma}
\def\ze{\zeta}
\def\io{\iota}
\def\bio{b^\iota}
\def\aio{a^\iota}
\def\twio{\tilde{w}^\iota}
\def\hwio{\hat{w}^\iota}
\def\gio{\g^\iota}
\def\Bio{B^\iota}

\def\del{\delta}
\def\pa{\partial}
\def\vp{\varphi}
\def\ve{\varepsilon}
\def\inf{\infty}

\def\vph{\varphi}
\def\vps{\varpsi}
\def\vPh{\varPhi}
\def\vep{\varepsilon}
\def\vpi{{\varpi}}
\def\vth{{\vartheta}}
\def\vsi{{\varsigma}}
\def\vrh{{\varrho}}

\def\bph{\bar{\phi}}
\def\bsi{\bar{\si}}
\def\bvp{\bar{\varphi}}

\newcommand{\bS}{{\mathbf S}}
\newcommand{\bH}{{\mathbf H}}
\newcommand{\bF}{{\mathbf F}}
\newcommand{\bE}{{\mathbf E}}

\def\tal{\tilde{\alpha}}
\def\tbe{\tilde{\beta}}
\def\tde{\tilde{\delta}}
\def\tpi{\tilde{\pi}}
\def\txi{\tilde{\xi}}
\def\tPi{\tilde{\Pi}}
\def\tPhi{\tilde{\Phi}}
\def\tV{\tilde{V}}
\def\tJ{\tilde{J}}
\def\tla{\tilde{\lambda}}
\def\tga{\tilde{\gamma}}
\def\tGa{\tilde{\Gamma}}
\def\tvs{\tilde{{\varsigma}}}
\def\tu{\tilde{u}}
\def\tU{\tilde{U}}
\def\tw{\widetilde w}
\def\tW{\widetilde W}
\def\tB{\tilde B}
\def\tv{\tilde v}
\def\tV{\tilde V}
\def\tz{\tilde z}
\def\tb{\tilde b}
\def\ta{\tilde a}
\def\tih{\tilde h}
\def\trh{\tilde {\rho}}
\def\tx{\tilde x}
\def\tf{\tilde f}
\def\tg{\tilde g}
\def\tG{\tilde G}
\def\tk{\tilde k}
\def\tl{\tilde l}
\def\tL{\tilde L}
\def\tD{\tilde D}
\def\tR{\tilde R}
\def\tP{\tilde P}
\def\tH{\tilde H}
\def\tp{\tilde p}

\def\hH{\hat{H}}
\def\hh{\hat{h}}
\def\hR{\hat{R}}
\def\hY{\hat{Y}}
\def\hX{\hat{X}}
\def\hP{\hat{P}}
\def\hT{\hat{T}}
\def\hV{\hat{V}}
\def\hG{\hat{G}}
\def\hF{\hat{F}}
\def\hw{\widehat{w}}
\def\hW{\widehat{W}}
\def\hu{\hat{u}}
\def\hs{\hat{s}}
\def\hv{\hat{v}}
\def\hb{\hat{b}}
\def\hB{\widehat{B}}
\def\hze{\hat{\zeta}}
\def\hsi{\hat{\sigma}}
\def\hrh{\hat{\rho}}
\def\hth{\hat{\theta}}
\def\hy{\hat{y}}
\def\hx{\hat{x}}
\def\hz{\hat{z}}
\def\hg{\hat{g}}
\def\he{\hat{e}}
\def\hE{\widehat{E}}

\def\B{\mathbf{B}}
\def\I{\mathbf{I}}
\def\P{\mathbf{P}}
\def\G{\mathbf{G}}
\def\S{\mathbf{S}}
\def\F{\mathbf{F}}
\def\one{\mathbf{1}}
\def\Sn{\mathbf{S}_n}
\def\0{\mathbf{0}}
\def\H{\mathbf{H}}
\def\V{\mathbf{V}}

\def\f{\mathcal{F}}
\def\çF{\mathcal{F}}
\def\o{\mathcal{O}}
\def\t{\mathcal{T}}
\def\r{\mathcal{R}}
\def\l{\mathcal{L}}
\def\m{\mathcal{M}}
\def\k{\mathcal{K}}
\def\n{\mathcal{N}}
\def\d{\mathcal{D}}
\def\p{\mathcal{P}}
\def\cP{\mathcal{P}}
\def\a{\mathcal{A}}
\def\h{\mathcal{H}}
\def\c{\mathcal{C}}
\def\y{\mathcal{Y}}
\def\e{\mathcal{E}}
\def\v{\mathcal{V}}
\def\z{\mathcal{Z}}
\def\x{\mathcal{X}}
\def\s{\mathcal{S}}
\def\g{\mathcal{G}}
\def\u{\mathcal{U}}
\def\w{\mathcal{W}}
\def\i{\mathcal{I}}
\def\j{\mathcal{J}}
\def\b{\mathcal{B}}

\def\lan{\langle}
\def\llb{(\!(}
\def\ran{\rangle}
\def\rrb{)\!)}
 \def\dim{{\hbox{\rm dim}}_{\mathbb C}\,}
\def\lng{\hbox{\rm{\tiny lng}}}
\def\sht{\hbox{\rm{\tiny sht}}}
\def\sph{\hbox{\rm{\tiny sph}}}
\def\inv{\hbox{\rm{\tiny inv}}}

\def\br#1{\langle #1 \rangle}

\def\rank{\hbox{rank}}
\def\gl{\mathfrak{gl}_N}

\newcommand{\Aut}{\operatorname{Aut}}
\newcommand{\Hom}{\operatorname{Hom}}
\newcommand{\End}{\operatorname{End}}
\newcommand{\Ind}{\operatorname{Ind}}
\newcommand{\ad}{\operatorname{ad}}
\newcommand{\pr}{\operatorname{pr}}
\newcommand{\aweyl}{\tilde{\mathbb S}_n}
\newcommand{\hec}{{\mathcal H}^t_n}
\newcommand{\Func}{{\mathcal F}({\mathbb C}^n,{\mathcal H}^t_n)}
\newcommand{\tr}{\operatorname{tr}}
\newcommand{\Out}{\operatorname{Out}}
\newcommand{\Rad}{\operatorname{Rad}}
\newcommand{\Spec}{\operatorname{Spec}}
\newcommand{\id}{\operatorname{id}}
\newcommand{\Int}{\operatorname{Int}}
\newcommand{\ct} {\operatorname{ct}}

\newcommand{\rat}{{\mathbb Q}}
\newcommand{\real}{{\mathbb R}}
\newcommand{\cplx}{{\mathbb C}}
\newcommand{\zint}{{\mathbb Z}}

\newcommand{\sq}{\phantom{1}\hfill$\qed$}
\newcommand{\Rea}{\Re}
\newcommand{\Ima}{\Im}

\newcommand{\st}{\bowtie}
\newcommand{\modd}{\mbox{\,mod\,}}
\newcommand{\lr}{\langle}
\newcommand{\rr}{\rangle}
\newcommand{\eps}{\varepsilon}
\newcommand{\phk}{\phi^{(k)}}
\newcommand{\psk}{\psi^{(k)}}
\newcommand{\Res}{\mbox{Res}\;}
\newcommand{\sgn}{\mbox{sgn}}
\newcommand{\mn} {\left\{ \begin{array}{c}m\\
n\end{array}\right\}}

\def\TT{\mathfrak{T}}
\def\JJ{\mathfrak{J}}
\def\HH{\mathfrak{H}}
\def\FF{\mathfrak{F}}
\def\GG{\mathfrak{G}}
\def\CC{\mathfrak{C}}
\def\LL{\mathfrak{L}}

\def\BB{\mathfrak{B}}
\def\AA{\mathfrak{A}}
\def\ZZ{\mathfrak{Z}}
\def\HH{\hbox{${\mathcal H}$\kern-5.2pt${\mathcal H}$}}
\def\tHH{\widetilde{\HH\ }}

\font\smm=msbm10 at 12pt
\def\symbol#1{\hbox{\smm #1}}
\def\lsmash{{\symbol n}}
\def\rsmash{{\symbol o}}
\def\#{\sharp}

\font\tenbf=cmbx10
\font\tenrm=cmr10
\font\tenit=cmti10
\font\ninebf=cmbx9
\font\ninerm=cmr9
\font\nineit=cmti9
\font\eightbf=cmbx8
\font\eightrm=cmr8
\font\eightit=cmti8
\font\sevenrm=cmr7
\font\sevenbf=cmbx7


\renewcommand{\natural}{\wr}

\numberwithin{equation}{section}

\section{Introduction}

This paper is a continuation of the part of \cite{C0}
devoted to {\em non-gatherable triangle triples},
NGT, in $\la$-sequences. The latter are the
sequences of positive roots associated with
reduced decompositions (words) in af\/f\/ine and non-af\/f\/ine
Weyl groups. We demonstrate that minimal NGT can be completely
described in the non-af\/f\/ine case; the answer appears especially
simple for the classical root systems and for $F_4$, $E_6$
(there are no NGT for $A_n$, $B_2$, $C_2$, $G_2$). As for $F_4$, $E_{6,7,8}$,
we reduced the problem to certain verif\/ications performed
by computer; it will be discussed in further works,
as well as af\/f\/ine generalizations.

The existence of NGT is a combinatorial obstacle for
using the technique of intertwiners (see, e.g.~\cite{C0})
for an {\em explicit} description of the irreducible representations
of the af\/f\/ine (and double af\/f\/ine) Hecke algebras, complementary to the
geometric theory of~\cite{KL1}. However, NGT are interesting
in their own right. Gathering together the triangle triples
using the Coxeter transformations seems an important question in the
theory of reduced decompositions of Weyl groups, which is far
from being simple.
More generally, assuming that
$\la(w)$ contains all positive roots of a root subsystem, can they
be gathered using the Coxeter transformations?

Let $R\in \R^n$ be a reduced irreducible root system or
its af\/f\/ine extension, $W$ the corresponding Weyl group.
Then the $\la$-{\em set} is def\/ined as
$\la(w)=R_+\cap w^{-1}(-R_+)$ for $w\in W$, where
$R_+$ is the set of positive roots in $R$. It is well-known
that $w$ is uniquely determined by
$\la(w)$; many properties of $w$ and its
reduced decompositions can be interpreted in terms
of this set. The $\la$-{\em sequence}
is the $\la$-set with the ordering of roots naturally
induced by the corresponding reduced decomposition.

The intrinsic description of such sets and sequences
is given in terms of the {\em triangle triples}
$\{\be,\ga=\al+\be,\al\}$. For instance, $\al,\be\in \la(w)$
$\Rightarrow$ $\al+\be\in \la(w)$ and the latter
root must appear between $\al$ and $\be$\, if this set
is treated as a sequence. This property is necessary but
not suf\/f\/icient; here and below see \cite{C0} for a
comprehensive discussion.

We need to know when a set of positive
roots of a rank two subsystem inside a
given {\em sequence}~$\la(w)$ can be {\em gathered} (made consecutive)
using the Coxeter transformations in $\la(w)$. It is natural
to allow the transformations only within
the minimal segment containing these roots.
This problem can be readily
reduced to considering the {\em triangle triples}
$\{\be,\ga=\al+\be,\al\}$ provided some
special conditions on the lengths.
The answer is always af\/f\/irmative
only for the root systems $A_n$, $B_2$, $C_2$, $G_2$ (and their
af\/f\/ine counterparts) or in the case when $|\al|\neq|\be|$;
otherwise NGT always exist.

For the root system $A_n$, gathering the triples is simple.
It readily results from the planar interpretation of
the reduced decompositions and the corresponding
$\la$-sequences in terms of $(n+1)$ lines in the
two-dimensional plane. This interpretation is essentially
equivalent to the classical {\em geometric} approach to the
reduced decompositions of $w\in W$ in terms of the lines
(or pseudo-lines) that go
from the main Weyl chamber to the chamber
corresponding to~$w$; see~\cite{Bo}.

The $A_n$-planar interpretation  was extended in \cite{Ch0}
to other {\em classical} root systems and $G_2$,
and then to their af\/f\/ine extensions in \cite{Ch5}.
It is given in terms of $n$ lines in $\R^2$ with ref\/lections
in one or two ``mirrors'' for $B_n$, $C_n$, $D_n$ (two mirrors are
needed in the af\/f\/ine case) or in terms of~$(n+1)$ lines
on the two-dimensional cylinder for the af\/f\/ine $A_n$.

We use the planar interpretation for the non-af\/f\/ine systems
$B$, $C$, $D$ to f\/ind {\em all} minimal non-gatherable triples,
{\em minimal NGT}, in these three cases.
No such interpretation is known for~$F_4$,~$E_{6,7,8}$, but
we managed to calculate all minimal NGT in these cases as well.
The af\/f\/ine root systems will be considered in the next paper.

Generally, the {\em admissibility} condition from \cite{C0} is
necessary and suf\/f\/icient for the triple to be {\em gatherable},
which is formulated in terms of subsystems of $R$ of types
$B_3$, $C_3$ or $D_4$.
We (re)establish this theorem in the non-af\/f\/ine case
in this paper and make the proof very constructive.
The proof presented in \cite{C0} was entirely algebraic, not quite
complete for the system $F_4$ and sketchy in the $D,E$-cases.

It is important to note that
the existence of NGT and other facts of similar nature are in
sharp contrast with the case of $A$. Generally, the theory
of root systems is uniform at level of generators and
relations of the corresponding Weyl or braid group;
however the root systems behave dif\/ferently when the
``relations for Coxeter relations'' are considered, i.e.,
at level of the second fundamental group.

Presumably, the phenomenon of NGT is one of the major
combinatorial obstacles for creating a universal
theory of AHA-DAHA ``highest vectors''
generalizing Zelevinsky's segments in the $A$-case
and based on the intertwining operators.
This technique was fully developed
only for af\/f\/ine and double af\/f\/ine Hecke algebras
of type $A_n$ and in some cases of small ranks.

The classif\/ication and explicit description
of {\em semisimple} representations of AHA and DAHA
is a natural application of this technique.
The fact that all triples are
gatherable in the case of~$A_n$ was important in~\cite{Ch1}  and
in further papers on the {\em quantum fusion procedure};
in type $A$,
AHA and DAHA are dual to quantum groups and quantum
toroidal algebras, generalizing
{\em affine Kac--Moody algebras}.
Extending the technique of intertwiners to other root systems
requires a~thorough analysis of NGT.

\section{Weyl groups}

Let $R=\{\al\}   \subset \R^n$ be a root system of type
$A,B,\dots,F,G$
with respect to a Euclidean form $(z,z')$ on $\R^n
\ni z,z'$,
$W$ the {\em Weyl group} \index{Weyl group $W$}
generated by the ref\/lections $s_\al$,
$R_{+}$ the set of positive  roots
corresponding to f\/ixed simple
roots $\al_1,\dots,\al_n$,
$\Ga$ the Dynkin diagram
with $\{\al_i,\, 1 \le i \le n\}$ as the vertices.
We will sometimes use
the dual roots (coroots) and the dual root system:
\[
R^\vee=\{\al^\vee =2\al/(\al,\al)\}.
\]


Let $\th\in R^\vee $ be the {\em maximal positive
root}, $\vth\in R^\vee $ the maximal {\em short} root. The
latter root is also the maximal positive coroot under the
normalization  $(\al,\al)=2$ for {\em short} roots.
Recall that
$1\ge(\th,\al^\vee)\ge 0$ for $\th\neq\al>0$. Similarly,
$1\ge(\vth,\al^\vee)\ge 0$ for $\vth\neq\al>0.$

Note that the sum of the long roots is always long,
the sum of two short roots can be a long root only
if they are orthogonal to each other.

The {\em length} of the reduced decomposition of $w\in W$ in
terms of the simple ref\/lections $s_i=s_{\al_i}$ $(1 \le i \le n)$
is denoted by $l(w)$. It can be also def\/ined as the
cardinality $|\la(w)|$
of the {\em $\la$-set} of $w$:
\begin{gather}\label{lasetdef}
 \la(w)\equal R_+\cap w^{-1}(R_-)=\{\al\in R_+,\
w(\al)\in R_-\}, \qquad
w\in W.
\end{gather}


The coincidence with the previous def\/inition
is based on the equivalence of the {\em length equality}
\begin{gather}\label{ltutwa}
 (a)\ \ l_\nu(wu)=
l_\nu(w)+l_\nu(u)
\qquad {\rm for}\quad w,u\in W
\end{gather}
and the {\em cocycle relation}
\begin{gather}
   (b)\ \ \la_\nu(wu) = \la_\nu(u) \cup
u^{-1}(\la_\nu(w)),
\label{ltutw}
\end{gather}
which, in its turn, is equivalent to
the {\em positivity condition}
\begin{gather}\label{ltutwc}
 (c)\ \  u^{-1}(\la_\nu(w))
\subset R_+.
\end{gather}
Applying (\ref{ltutw}) to the reduced decomposition
$w=s_{i_l}\cdots s_{i_2}s_{i_1}$:
\begin{gather}
\la(w) = \{ \al^l=w^{-1}s_{i_l}(\al_{i_l}), \;
\ldots, \;\al^3=s_{i_1}s_{i_2}(\al_{i_3}), \; \al^2=s_{i_1}(\al_{i_2}),\;  \al^1=\al_{i_1}   \}.
\label{tal}
\end{gather}
This relation demonstrates directly that the cardinality
$l$ of the set $\la(w)$ equals $l(w)$.
Cf.~\cite[Section~4.5]{Hu}.
We also note that
$\la_\nu(w^{-1}) = -w(\la_\nu(w))
$.

It is worth mentioning that a counterpart of the
$\la$-set can be introduced for reduced decomposition
$w=s_{i_l}\cdots s_{i_2}s_{i_1}$
in arbitrary Coxeter groups. Following \cite[Ch. IV, 1.4, Lemma 2]{Bo}
 one can def\/ine
\begin{gather}
\La(w) = \{t_l=w^{-1}s_{i_l}(s_{i_l}),\;
\ldots,\; t_3=s_{i_1}s_{i_2}(s_{i_3}),
t_2=s_{i_1}(s_{i_2}),\; t_1=s_{i_1}  \},
\label{talc}
\end{gather}
where the action is by conjugation; $\La(w)\subset W$.

The $t$-elements are pairwise dif\/ferent if and
only if the decomposition is reduced (a simple
straight calculation; see \cite{Bo});
then this set does not depend
on the choice of the reduced decomposition. It readily
gives a proof of formula~(\ref{tal}) by induction
and establishes the equivalence
of~(a), (b) and (c).

Using the root system  dramatically simplif\/ies
theoretical and practical (via computers) ana\-ly\-sis of
the reduced decompositions and makes the
crystallographical case signif\/icantly simpler
than the case of abstract Coxeter groups.
The positivity of roots,
the alternative def\/inition of the $\la$-sets from
(\ref{lasetdef}) and, more specif\/ically,
property (c) are of obvious importance.
These features are (generally)
missing in the theory of abstract Coxeter groups,
though the $\La$-sets from~(\ref{talc})
can be of course used for various questions.

The sets  $\la(w)$ can be treated naturally as
a {\em sequence}; the roots in (\ref{tal}) are ordered
naturally depending on the particular choice
of a reduced decomposition.
We will mainly treat $\la(w)$ as
sequences in this paper,
called $\la$-{\em sequences}.

Note that relation (\ref{ltutwc}) readily gives that
an arbitrary simple root $\al_i\in\la(w)$ can be made
the f\/irst in a certain $\la$-sequence. More generally:
\begin{gather}
 \la_\nu(w)
=  \{\al>0, \; l_\nu( w s_{\al}) \le l_\nu(w) \};
\label{xlambda1}
\end{gather}
see \cite{Bo}
and \cite[Section~4.6, Exchange Condition]{Hu}. This property is
closely related to the formula:
\begin{gather} \label{talinla}
 \al\in \la(w) \Leftrightarrow
\la(s_{\al}) =  \{\be,\;-s_{\al}(\be)\,\mid\,
s_{\al}(\be)\in -R_+ ,\; \be\in\la(w)\}.
\end{gather}

\section{Coxeter transformations}

We will prepare some tools for studying transformations of
the reduced decompositions. The elementary ones are the
{\em Coxeter transformations} that are substitutions
$(\cdots s_is_js_i)\mapsto (\cdots s_js_is_j)$ in
reduced decompositions of the elements $w\in W$; the number
of $s$-factors is  2, 3, 4, 6  as $\al_i$ and $\al_j$ are connected
by $m_{ij}=0,1,2,3$ laces in the af\/f\/ine or non-af\/f\/ine
Dynkin diagram. They induce
{\em reversing the order} in the corresponding segments
(with  2, 3, 4, 6 roots)
inside the sequence $\la(w)$. The corresponding roots form
a set identif\/ied with the set of positive roots
of type $A_1\times A_1$, $A_2$, $B_2$, $G_2$ respectively.

The theorem below is essentially from \cite{C0}; it
has application to the decomposition of the
polynomial representation of DAHA, the classif\/ication
of semisimple representations of AHA, DAHA
and to similar questions. We think that it clarif\/ies
why dealing with the intertwining
operators for arbitrary root systems is signif\/icantly more
dif\/f\/icult than in the $A_n$-theory (where much is known).

Given a reduced decomposition of $w\in W$,
let us assume that $\al+\be=\ga$ for the roots
$\ldots,\be,\ldots,\ga,\ldots,\al\ldots$
in $\la(w)$ ($\al$ appears the f\/irst), where
only the following combinations of their lengths
are allowed in the $B$, $C$, $F$ cases
\begin{gather}
\hbox{lng}+\hbox{lng}=\hbox{lng}\ \,
 (B,F_4)\qquad \hbox{or}\qquad
\hbox{sht}+\hbox{sht}=\hbox{sht} \ \,
(C,F_4).
\label{shtshtsht}
\end{gather}
Since we will use the Coxeter transformations
only inside the
segment $[\be,\al]\subset\la(w)$ between $\be$ and
$\al$, it suf\/f\/ices to assume that  $\al$ is a simple root.
Also, let us exclude $A_n$, $B_2$, $C_2$, $G_2$ from the consideration
(in these cases all triangle triples, if any,
are gatherable).

\begin{theorem}\label{RANKTWO}
(i) For the root systems of type $B_n$, $C_n$, $F_4$,
the roots $\be$, $\ga$, $\al$ are non-gatherable
(cannot be made consecutive
using the Coxeter transformations) inside the segment
$[\be,\al]\subset\la(w)$ if and only if a root
subsystem $R^3\subset R$ of type
$B_3$ or $C_3$ ($m=1,2$) exists such that
\begin{align}\label{rankthrees}
&\be=\ep_1+\ep_3,\, \al=\ep_2-\ep_3,\
\ep_1-\ep_2\,\not\in\, [\,\be,\al\,]\,
\not\ni\, m\ep_3,
\end{align}
where the roots $\ep_1, \ep_2,\ep_3 \in R^3$ are from the
$B_3$, $C_3$ tables of {\rm \cite{Bo}}. Equivalently,
the sequence $[\be,\al]\cap R^3_+$ (with the
natural ordering) must be
\begin{align}
&\{\,\ep_1+\ep_3,\, m\ep_1,\, \ep_2+\ep_3,\, \ga=\ep_1+\ep_2,\,
\ep_1-\ep_3,\,
m\ep_2,\, \ep_2-\ep_3\,\} \label{rankthree}
\end{align}
up to Coxeter transformations in $R^3$ and changing the order
of all roots in \eqref{rankthree} to the opposite.
This sequence is $\la^3(s_{\ga})$ in $R_+^3$ for
the maximal root $\ga=\th^3$ for $B_3$ and for
the maximal short root  $\ga=\vth^3$ for $C_3$.

(ii) For the root system
$R$ of type $D_{n\ge 4}\,$ or for
$E_{6,7,8}$, $\{\be,\ga,\al\}$ is a non-gatherable triple
if and only if
a root subsystem $R^4\subset R$ of type
$D_4$ can be found such that
\begin{gather}\label{rankthreesd}
\be=\ep_1+\ep_3,\qquad \ga = \ep_1+\ep_2,\qquad
\al = \ep_2-\ep_3,
\\
\{ \ep_1-\ep_2,\,\ep_3-\ep_4,\,
\ep_3+\ep_4 \}\cap [ \be,\al ]  =\varnothing,\notag
\end{gather}
where $\ep_1$, $\ep_2$, $\ep_3$, $\ep_4$ are from the $D_4$-table of
{\rm \cite{Bo}}. Equivalently,
the sequence $[ \be,\al ]\cap R^4_+$  must be
\begin{gather}
\{ \be=\ep_1+\ep_3,\, \ep_1-\ep_4,\, \ep_1+\ep_4,\,
\ep_1-\ep_3, \ga=\ep_1+\ep_2,\notag\\
 \phantom{\{ \be=\ep_1+\ep_3,\,}{} \ep_2+\ep_3,\, \ep_2+\ep_4,\, \ep_2-\ep_4,\,
\al=\ep_2-\ep_3 \}\label{rankthreed}
\end{gather}
up to Coxeter transformations in $R^4$. Equivalently,
$[\,\be,\al\,]\cap R^4_+$
is the $\la$-set of $s_{\th^4}$ in $R_+^4$ for
the maximal root $\th^4$.
\end{theorem}

We will (re)prove this theorem (later) by listing all
{\em minimal non-gatherable triples}. Our approach
is signif\/icantly more explicit than
that from \cite{C0}, although Theorem~\ref{RANKTWO} does not require the minimality
condition and therefore is of more general nature.
The af\/f\/ine root systems will be considered elsewhere.

To begin with,
the following are the lists of the non-af\/f\/ine roots
$\ga>0$ such that the endpoints $\be=\ga-\al_j$, $\al=\al_j$ of
$\la(s_{\ga})$ are {\em non-movable} under
the Coxeter transformations within $\la(s_{\al})$ and
$\{\be,\ga,\al\}$ form an $A_2$-triple in the cases of
$F_4$, $B$, $C$; $A_2$-triples are those
subject to  $|\al_j|=|\ga|=|\be|$.
See \cite{C0} and also (\ref{talinla}).
The maximal long root $\th$ (for $B_{n\ge 3}$, $F_4$ and
in the simply-laced case) and maximal short root $\vth$
(for $C_{n\ge 3}$, $F_4$) are
examples of such $\ga$ (but there are many more).

The bar will show the position of the corresponding $\al_j$ in the
Dynkin diagram. We will omit
the cases of $E_{7,8}$; there are
$7$ such $\ga$ for $E_7$ and $22$ for $E_8$.

{\bf The case of {\mathversion{bold}$E_6$}.}
The roots $\ga\in R_+$ such that $\la(s_\ga)$ has
non-movable endpoints are:
\begin{align}\label{e6sing}
01&\overline{2}10,& 1\overline{2}&210,& 01&2\overline{2}1,&
12&\overline{3}21,& 12&321.         & \\
  &1              &              &1   &   &1              &
  &1              &   &\overline{2} & \notag
\end{align}
The corresponding triple $\{\be=\ga-\al_j,\ga,\al_j\}$
is a minimal non-gatherable triple inside $\la(s_{\ga})$.

{\bf The case of {\mathversion{bold}$F_4$}.}
The roots $\ga\in R_+$ with non-movable endpoints
of $\la(s_\ga)$ and
subject to $|\al_j|=|\ga|$ are:
\begin{align}\label{f4sing}
&01\overline{2}1,& &1\overline{2}20,& &12\overline{3}1,&
&123\overline{2},& &1\overline{3}42,& &\overline{2}342.&
\end{align}

{\bf The case of {\mathversion{bold}$B$, $C$, $D_n$}.}
Given $\al_j=\ep_j-\ep_{j+1}$, the corresponding
root $\ga$ (it is unique)
equals  $\ep_{j-1}+\ep_{j}$ for $j=2,\ldots,n-1$  provided
that $n\ge 3$ and $j<n-1\ge 3$ for $D_n$.
The notation is from~\cite{Bo}.

\section{Minimal NGT}

We are now in a position to formulate the main result
of the paper, that is a description of all minimal
non-gatherable triples, NGT, for the non-af\/f\/ine root
systems. It provides a direct justif\/ication
of Theorem \ref{RANKTWO} in the non-af\/f\/ine case.
We will omit the lists in the cases $E_{7,8}$
in this paper (there are no NGT for $A_n$, $B_2$, $C_2$, $G_2$).

We say that $w\in W$ {\em represents a minimal NGT}, if
$\{\be,\ga=\al+\be,\al\}\subset \la(w)$,
$\al$, $\be$ are correspondingly the beginning and
the end of the sequence $\la(w)$ and these roots are non-movable
inside $\la(w)$ using the Coxeter transformations.
Recall that the condition from (\ref{shtshtsht}) is imposed.

\begin{theorem}\label{THMBCFE}
The lists of elements $w\in W$ representing minimal
NGT are as follows.

(i) For $B_n$, $C_n$, $D_n$ an arbitrary
$\ga=\ep_i+\ep_j, i<j$, subject to $j<n$ for $B,C$
and  $j<n-1$ for~$D$ can be taken; the corresponding
simple $\al$ will be $\ep_{j-1}-\ep_{j}$ in the notation
from {\rm \cite{Bo}}. The element $w$ is the product of
reflections corresponding to the ``telescoping'' sequence
$\ep_i+\ep_j, \ep_{i+1}+\ep_{j-1},\ldots$ ending with
$\ep_{k}+\ep_{k+1}$ unless the last root of this sequence
is in the form $\ep_{k-1}+\ep_{k+1}$. In the latter case,
the roots $\ep_k$ or  $2\ep_k$ must be added to this sequence
for $B_n$ or $C_n$, and the pair of roots $\ep_k+\ep_n$, $\ep_k-\ep_n$
must be added for $D_n$.

Such $w$ is determined
uniquely by $\{i,j\}$ and is a product of reflections for
pairwise orthogonal roots; in particular, $w^2={\rm id}$.
One of these roots must be $\ga$ (cf.\ the description
of $w=s_{\ga}$ given above).

(ii) In the case of $F_4$, such $w$ are products of
pairwise commutative reflection as well, but $w$ is
not uniquely determined by the triple. Omitting two
$w$ that come from $B_3$ and $C_3$ naturally embedded
into $F_4$, the remaining eight are as follows: there are
four that are simply reflections with respect to the roots
$1231$, $1342$, $1232$, and $2342$ from \eqref{f4sing};
the remaining four are each the
product of reflections of two orthogonal roots,
$s_{1342}\,s_{1110}$,
$s_{0122}\,s_{1231}$, $s_{1222}\,s_{1231}$, and
$s_{1342}\,s_{1121}$.
Note that since $F_4$ is self dual, the dual of every word on
this list is also on this list.

(iii) In the case of $E_6$ we will omit $5$ elements
 coming from the two natural  $D_5$-subsystems of $E_6$.
($3+3$ minus the one for $D_4$ they have in common; see~(i)).
Of the remaining ten, two are reflections of the roots
$\begin{array}{c}
12321 \\
1
\end{array}$ and $\begin{array}{c}
12321 \\
2
\end{array}$ from \eqref{e6sing}.

Two more can be written as the composition
of three pairwise orthogonal reflections of roots:
\[
s\left(\begin{array}{c}
12321 \\
1
\end{array}\right)\,s\left(\begin{array}{c}
01100 \\
1
\end{array}\right)\,s\left(\begin{array}{c}
00110 \\
1
\end{array}\right), \qquad  s\left(\begin{array}{c}
01221 \\
1
\end{array}\right)\,s\left(\begin{array}{c}
12210 \\
1
\end{array}\right)\,s\left(\begin{array}{c}
11211 \\
1
\end{array}\right).
\]
The final six can not be written
as products of orthogonal reflections. Written as
products of reflections of simple roots they are as follows:
\begin{alignat*}{3}
&21324354632143254363243, \qquad &&  32143263214325436321432, & \\
&32435463214325436324354, \qquad &&  43215432632143254363243,& \\
&2132435463214325436324354, \qquad &&  4321543263214325436321432,&
\end{alignat*}
where we abbreviate $s_{\alpha_i}$  as  $i$
for $1 \leq i \leq 6$; for instance, the first word is
$ s_2s_1\cdots s_4s_3$.
\end{theorem}

\begin{proof}
We will use the planar interpretation of the reduced decompositions
from \cite{Ch0} for~$B$,~$C$,~$D$.
An arbitrary element $w\in W$ can be represented by a conf\/iguration
of $n$ lines in the plane with a~possible ref\/lection in the
$x$-axis. The initial numeration of the lines is from the top
to the bottom (the right column of the line numbers in
the f\/igures below). The intersections and ref\/lections
are supposed to have pairwise dif\/ferent $x$-projections;
simultaneous multiple {\em events}
(like double intersections) are not allowed.

Given an intersection, we plot the vertical line
through it and count the lines (locally) from the top;
the intersection of the (neighboring)
lines $k$, $k+1$ is interpreted as $s_k$.
The angle between these lines gives
the corresponding root in $\la(w)$. Namely, it
is $\ep_i\pm\ep_j$ for the initial (global)
numbers of the lines and their initial
angles $\ep_i$, $\ep_j$ with the $x$-axis; the sign depends
on the number of ref\/lections before the intersection.
See Fig.~\ref{bngt}, where we put $i$ instead of $\ep_i$ in
the angles/roots and instead of $s_i$ in the decomposition.

\begin{figure}[th]
\centerline{\includegraphics[scale=0.5]{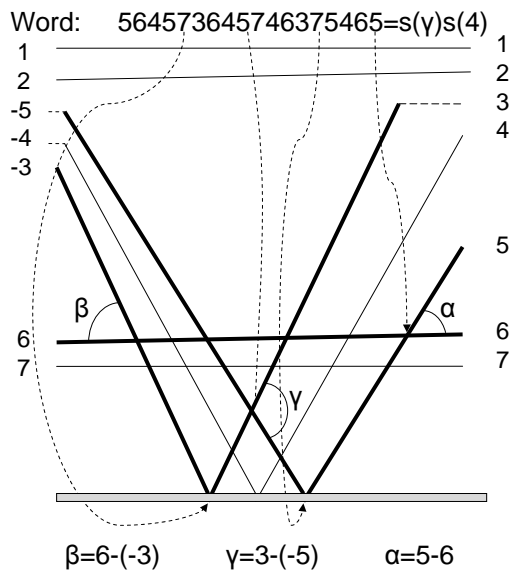}}

\caption{Typical minimal NGT for $B_7$.}\label{bngt}
\end{figure}

The angle is taken $\ep_i$ or $2\ep_i$ for the ref\/lection
in the $x$-axis in the cases of  $B$ or $C$; the corresponding
event is interpreted as $s_n$ in the Weyl group.

Treating the ref\/lections is a bit more involved in the
$D$-case. The combination of the
ref\/lection, then the $\{n-1,n\}$-intersection
(the numbers of lines are local), and then the ref\/lection
again is interpreted as $s_n$ for $D_n$.
The corresponding root from $\la(w)$
is the middle angle in this event, which
will be called {\em $V\times V$-shape}. These
events are encircled in Fig.~\ref{dngt}; they look like~{$\bigvee$}\kern -6.5pt {$\bigvee$}. Their angles
are $5+6$ ($\ep_5+\ep_6$, to be exact) and $3+4$
correspondingly (from right to left).

\begin{figure}[t]
\centerline{\includegraphics[scale=0.5]{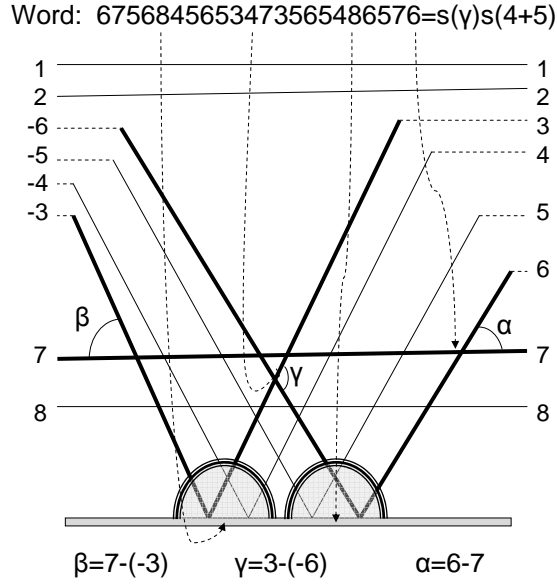}}
\caption{Typical minimal NGT for $D_8$.}
\label{dngt}
\end{figure}

This construction is suf\/f\/icient for constructing reduced
decompositions for arbitrary
conf\/i\-gu\-ra\-tions with the even number of ref\/lections
in the $x$-axis.
Indeed, proper moving the lines upward
makes any given conf\/iguration a sequence of the simple crosses
of lines and the $V\times V$-shapes.
However, the geometric interpretation of the Coxeter relation
$s_{n-2}s_n s_{n-2}=s_ns_{n-2}s_n$ requires an extra
{\em $V+V$-operation},
that is breaking a given line twice
and adding two ref\/lections,
as shown in Figs.~\ref{drelat} and~\ref{dngt1},
followed by creating the $V\times V$-shapes.
Symbolically, it can be represented by~{$\bigvee$}\kern -1pt{$\bigvee$} (line~4 in Fig.~\ref{drelat}).
More formally,
\begin{enumerate}\itemsep=0pt
\item[1)] given a line, this transformation must not increase the total
number of its intersections with the other lines;

\item[2)] two ref\/lections must exist in a given conf\/iguration
neighboring to the (new) ref\/lections from $V+V$;

\item[3)] the pairs of neighboring ref\/lections
from (2) have to be arranged into two $V\times V$-shapes.
\end{enumerate}

Performing one such $V+V$ or multiple operations of this type and
moving the lines if necessary, the f\/inal conf\/iguration can be
represented in terms of (simple) intersections and
$V\times V$-shapes, provided that the number of initial ref\/lections
is even. Fig.~\ref{drelat} gives the simplest minimal NGT
represented with and without transforming
line~4. Here one avoids breaking line~4 and
adding the $V+V$-shape to
this line by moving it up (the second picture).
Disregarding line~3, the f\/igure
represents the Coxeter relation $s_ns_{n-2}s_n=s_{n-2}s_ns_{n-2}$.

\begin{figure}[t]
\centerline{\includegraphics[scale=0.45]{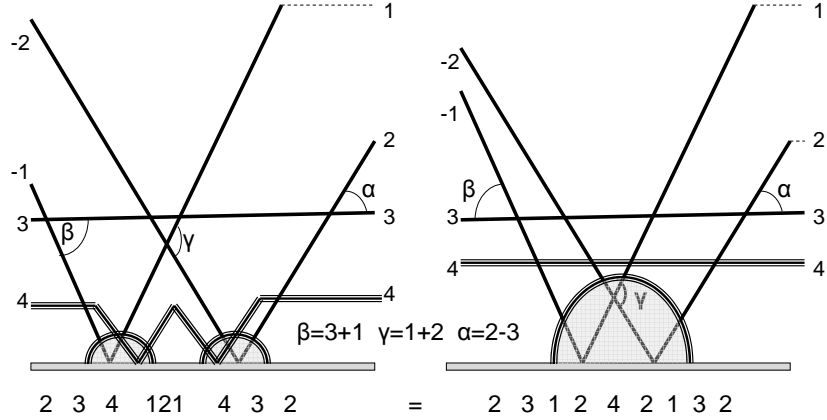}}
\caption{Transforming the line for $D_4$.}\label{drelat}
\end{figure}

The claim is that
{\em the resulting products of simple reflections (the words
in the figures) are always reduced; the angles
give the corresponding $\la$-sequence}. Recall,
that the ordering of the $x$-projections is from right to left
as well as that of the simple ref\/lections and angles.

It is natural to change
the global numbers of the lines from $i$ to $-i$ upon the
ref\/lections. Then the resulting column of global
line numbers (on the left)
gives the standard one-line presentation of the corresponding
$w$. For the $D$-system, the $V+V$ operation
does not change the global numbers at the ends,
since we changes the sign {\em two times} at the additional
ref\/lection points. Note that, technically, we do not change
the line number (the ``global angle'' assigned to this line)
at the beginning and at the end
of the additional $V+V$-shape;
these are ``no-events'', though, geometrically, the
angle of this line is changed at these points.

It is worth mentioning that not all reduced decompositions
of $w\in W$ can be obtained using the lines only; generally,
one should consider ``pseudo-lines'', where the assignment
of the ``global angle'' to the line is combinatorial rather
than geometric. Fig.~\ref{drelat} is a good demonstration
of this phenomenon (the counterexamples exist even for
$A_n$ with suf\/f\/iciently big $n$).

Using the planar interpretation,
the proof of (i) goes as follows.

First of all, any (triangle) triple for  $A_n$
can be readily made consecutive, corresponding
to a~``minimal'' triangle, using proper
parallel translations of the lines. The same argument
shows that the root $\ga$ cannot be $\ep_i-\ep_j$ for $B$, $C$, $D$
in the notation from \cite{Bo}. Otherwise, we can make the
corresponding triangle ``minimal'' as for $A_n$.

We will consider the $B, C$-case only; the root system $D$
is completely analogous. Given $\ga=\ep_i+\ep_j$,
there are three groups of the
(initial) lines:
\begin{enumerate}\itemsep=0pt
\item[(a)] beyond $i$ (lines  1, 2  in Fig.~\ref{bngt}),

\item[(b)] between $i$ and $j$ (lines  3, 4, 5  there) and

\item[(c)] the lines strictly below $j$ (namely, lines  6, 7).
\end{enumerate}

The lines from group (a) do not intersect each other.
Otherwise, we can {\em move} the intersection point to the right
using the parallel translations (as in the $A$-case) and
make the corresponding root the f\/irst in $\la(w)$,
which is impossible since $\al$ is the only such root.
Also the lines from (a) cannot intersect the
lines from group~(b). If such an intersection occurs then
we can move it to the right or to the left till the
very f\/irst or very last position in a reduced decomposition
of $w$, which contradicts to the minimality of the NGT
under consideration.

Similarly, the lines from group (c) cannot intersect each
other. They also do not intersect the lines from group (a).
Moreover, the intersections inside the group (b) can occur
only due to the ref\/lections (i.e., all pairs of lines in
this group must intersect exactly once).
If, say, lines~4 and~5 intersect before the ref\/lection
point of line $5$ or after the ref\/lection point of line~4
(see the $B$-f\/igure), then we can move this intersection to the
right or to the left all the way. Finally, the group (c)
is always nonempty (i.e., $j<n$); otherwise $s_n$ can be made
the f\/irst in a certain reduced decomposition of $w$.

It implies that the simple root $\al$ (from the triple
under consideration) can be only $\ep_j-\ep_{j+1}$,
i.e., the one corresponding to the intersection of
lines $j$ and $j+1$.
Indeed, the other possibility for~$\al$,
the simple root $\ep_i-\ep_{i+1}$, would
contradict the minimality of NGT. Respectively,
$\be=\ep_{i}+\ep_{j+1}$.

Summarizing, $w$ can be only of the type shown
in Fig.~\ref{bngt}. Geometrically, it is obvious that
$\{\be,\ga,\al\}$ there (and in general)
is a minimal non-gatherable triple. Indeed, the endpoints
are non-movable using the Coxeter transformations, which
correspond to moving (maybe with changing the angles)
the lines without altering their initial and f\/inal orderings,
i.e., the right and left columns of numbers.

The same reasoning gives that $j<n-1$ and that
minimal NGT can be only as shown in Fig.~\ref{dngt}
in the $D$-case. This concludes (i).

The lists (ii), (iii) are reduced to certain direct computer
calculations to be discussed elsewhere (including the complete
lists for $E_{7,8}$). We note
that f\/inding all $w$ representing minimal NGT
for $F_4$, $E_6$ is a relatively simple problem for
products of pairwise commutative ref\/lections (it is
not always true for $E_{6,7,8}$). It is a straightforward
generalizations of the description of the pure $w=s_\ga$
representing minimal NGT we gave above. One
of these ref\/lections must be $s_{\ga}$ for $\ga$ from
the triple; it simplif\/ies the consideration.
\end{proof}

\section[The existence of $R^{3,4}$]{The existence of $\boldsymbol{R^{3,4}}$}

The Theorem~\ref{RANKTWO} guarantees the existence
of the subsystems $B_3$ or $C_3$ and $D_4$,
ensuring that the corresponding (admissible)
triple is non-gatherable.
Recall that the intersection of the (positive roots of)
these subsystems with $\la(w)$ containing such
triple, must contain $7$ but not $9$ (the total)
roots for $B_3$, $C_3$ and $9$ but not $12$ (the total) roots
in the case of $D_4$.
We will call such $7$-subsets or $9$-subsets in $\la(w)$
{\em blocks for NGT}, respectively,
$B_3$-blocks, $C_3$-blocks, $D_4$-blocks.

The blocks can be naturally seen geometrically
in the cases $B_n$, $C_n$ and $D_n$.
Indeed, if one considers only
bold lines in Fig.~\ref{bngt}, then it readily leads
to the desired $R^3$ in the cases $B_n$, $C_n$. The intersection of the
$\la$-sequence with this $R^3$ will contain exactly $7$ roots
(from possible $9$), i.e., form a block, an obstacle for gathering
the corresponding triple.

\begin{figure}[t]
\centerline{\includegraphics[scale=0.5]{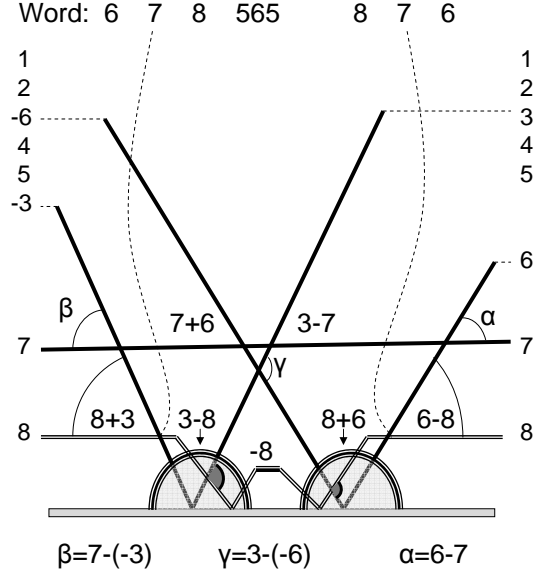}}
\caption{Finding $D_4$ in NGT for $D_8$.}
\label{dngt1}
\end{figure}

For f\/inding a root subsystem $R^4$ of
type $D_4$ in Fig.~\ref{dngt}, lines 3, 6, 7, 8 must
be taken. Line~8 must be moved up to make it
beyond $\ga$ (but below line~7) or transformed by adding
the $V+V$-shape. It is shown in Fig.~\ref{dngt1}.
If there are
several lines like 8 ``almost parallel'' to each other,
then either one can be taken to construct a $D_4$-block.

The number of the roots (only positive ones matter)
in the intersection of such $D_4$-subsystem
with the $\la$-sequence is always exactly $9$
(from possible $12$), i.e., this intersection
is a $D_4$-{\em block}.

We conclude that
{\em the $R^3$-subsystem and the corresponding block
is unique for a given non-gatherable
triple in types $B_n$, $C_n$}. In the case of $D_n$,
{\em a $D_4$-block always exists for NGT},  but the $R^4$-subsystem
is generally not unique. It proves
Theorem \ref{RANKTWO} for the classical roots systems
and makes explicit the constructions of {\em blocks}.

{\bf The case of {\mathversion{bold}$F_4$}.}
Given a word $w$ containing a minimal NGT where
all three roots are {\em long}, Theorem \ref{RANKTWO} states that
there exist three roots, $\alpha_1$, $\alpha_2$, $\alpha_3$,
with the following properties:

\begin{enumerate}\itemsep=0pt
    \item The roots $\alpha_1$, $\alpha_2$, $\alpha_3$
satisfy the properties of the set of simple roots of a $B_3$
root system using notation from \cite{Bo}.
    \item The NGT in $\la(w)$ can be written
$\be=\alpha_1 + \alpha_2 + 2\alpha_3$,
$\ga=\alpha_1 + 2\alpha_2 + 2\alpha_3$, $\al=\alpha_2$.
    \item The seven roots $\alpha_2$, $\alpha_1 + \alpha_2$,
$\alpha_2 + \alpha_3$, $\alpha_1 + \alpha_2 + \alpha_3$,
$\alpha_2 + 2\alpha_3$, $\alpha_1 + \alpha_2 + 2\alpha_3$, and
$\alpha_1 + 2\alpha_2 + 2\alpha_3$ form the intersection of
$\la(w)$ and the subsystem
$R^3\subset R$ generated by $\al_1$, $\al_2$, $\al_3$.
\end{enumerate}

Note that $\alpha_2$ is always a simple root of $F_4$
but~$\alpha_1$ and~$\alpha_3$ need not be.
The two roots~$\alpha_1$ and~$\alpha_3$ will not be
elements of $\lambda(w)$. Together with the seven roots above
these are all nine positive roots in the $R^3$, i.e.,
we come to the condition
{\em seven but not nine} given in~(\ref{rankthree}).

Following Theorem~\ref{THMBCFE}(ii),
we will explicitly demonstrate that Theorem \ref{RANKTWO}
holds for minimal NGT and give an appropriate choice
of $\alpha_1$, $\alpha_2$, $\alpha_3$ (there is often more than
one {\em block}). Consi\-de\-ring only {\em minimal} NGT is
obviously suf\/f\/icient to check  Theorem~\ref{RANKTWO}.

We begin with the following explicit example.
Let $w=s_{1342}s_{1110}$; the {\em shortlex form}
of $w$ is $2132132432132432=s_2s_1\cdots s_4s_3s_2$.
Here the ordering is {\em lexicographical} from left to right
(but we apply the corresponding ref\/lections from right to
left).

This yields:
\begin{gather*}
                            \lambda(w) =
 \{1242, 1120, 1232, 2342, 1222, 1110, 1100, 1231, \\
\phantom{\lambda(w) = \{}{}1221, 1342, 1220, 0121, 0120, 0111, 0110, 0100\}.
                 \end{gather*}
The NGT is $\{1242, 1342, 0100\}$.
If we choose \
$\alpha_1 = 1122$,  $\alpha_2 = 0100$,
$\alpha_3 = 0010$,  then
$\{\alpha_1 + \alpha_2 + 2\alpha_3,
\alpha_1 + 2\alpha_2 + 2\alpha_3,
\alpha_2\} = \{1242, 1342, 0100\}$ (the NGT), and
$\lambda(w)$ contains the $B_3$-block:
\begin{gather*}
\{ \alpha_2 = 0100, \alpha_1 + \alpha_2 = 1222,\,
\alpha_2 + \alpha_3 = 0110, \, \alpha_1 + \alpha_2 + \alpha_3 = 1232,\\
\phantom{\{}{} \alpha_2 + 2\alpha_3 = 0120,\,
\alpha_1 + \alpha_2 + 2\alpha_3 = 1242,\,
\alpha_1 + 2\alpha_2 + 2\alpha_3 =1342\}.
\end{gather*}
Finally, $\lambda(w)$ does not contain either $\alpha_1 = 1122$ or
$\alpha_3 = 0010$. So the seven but not nine condition is satisf\/ied.

\begin{table}[t]
\caption{$F_4$ Min-NGT's: $B_3$-blocks.}
\vspace{1mm}
\centering
{\small \begin{tabular}{|c|c|c|} \hline
$w$                                 & NGT
      & $\alpha_1$, $\alpha_2$, $\alpha_3$ \\ \hline
$s_{1220}$
& $1120, 1220, 0100$    & $1000, 0100, 0010$ \\ \hline
$s_{1342}$
& $1242, 1342, 0100$    & $1120, 0100, 0011$ \\ \hline
$s_{2342}$
& $1342, 2342, 1000$    & $0100, 1000, 0121$ \\ \hline
$s_{1342}s_{1110}$
& $1242, 1342, 0100$    & $1122, 0100, 0010$ \\ \hline
$s_{1222}s_{1231}$
& $1122, 1222, 0100$    & $1000, 0100, 0011$ \\ \hline
\end{tabular}}
\label{F4LongTable}
\end{table}

Table \ref{F4LongTable} shows each of the Min-NGT words in
$F_4$ where the NGT is made up of {\em long} roots and,
correspondingly, the {\em block} must be of type $B_3$.
Also included are the NGT and a choice of~$\alpha_1$, $\alpha_2$, $\alpha_3$ that determine an appropriate
$B_3$-block. The word used in the example above is also
included.

{\bf Short NGT for {\mathversion{bold}$F_4$}.}
Similarly, if the roots from
Min-NGT are all {\em short}, then Theorem \ref{RANKTWO} in the
case of $F_4$
tells us that there exist three roots,
$\alpha_1$, $\alpha_2$, $\alpha_3$, with the following properties:

\begin{enumerate}\itemsep=0pt
    \item The roots $\alpha_1$, $\alpha_2$, $\alpha_3$ satisfy the
properties of the set of simple roots of a $C_3$ root
system using notation from \cite{Bo}.
    \item The NGT in $\la(w)$ is written
$\alpha_1 + \alpha_2 + \alpha_3$, $\alpha_1 + 2\alpha_2 + \alpha_3$,
$\alpha_2$.
    \item The intersection $\la(w)\cap R^3$ is formed by
$\alpha_2$, $\alpha_1 + \alpha_2$,
$\alpha_2 + \alpha_3$, $\alpha_1 + \alpha_2 + \alpha_3$,
$2\alpha_2 + \alpha_3$, $\alpha_1 + 2\alpha_2 + \alpha_3$, and
$2\alpha_1 + 2\alpha_2 + \alpha_3$
for $R^3\subset R$
generated by $\alpha_1$, $\alpha_2$, $\alpha_3$.
\end{enumerate}

Here $\alpha_2$ is always a simple root of $F_4$ but $\alpha_1$
and $\alpha_3$ need not be.
The two roots $\alpha_1$ and $\alpha_3$ will not be
elements of $\lambda(w)$. Together with the seven roots above these
are all nine positive roots of $R^3$.
This condition {\em seven but not nine}
from (\ref{rankthree}) is satisf\/ied.

Table~\ref{F4ShortTable} shows each of the Min-NGT words in~$F_4$ where the NGT is made up of short roots.
Also included are the roots that make up the NGT and a choice of
$\alpha_1$, $\alpha_2$, $\alpha_3$ that determines an appropriate
$C_3$ subsystem.

\begin{table}[t]
\caption{$F_4$ Min-NGT: $C_3$-blocks.}
\vspace{1mm}

\centering
{\small \begin{tabular}{|c|c|c|} \hline
$w$                                 & NGT
       & $\alpha_1$, $\alpha_2$, $\alpha_3$ \\ \hline
$s_{0121}$
& $0111, 0121, 0010$    & $0001, 0010, 0100$ \\ \hline
$s_{1231}$
& $1221, 1231, 0010$    & $0111, 0010, 1100$ \\ \hline
$s_{1232}$
& $1231, 1232, 0001$    & $0010, 0001, 1220$ \\ \hline
$s_{1231}s_{0122}$
& $1221, 1231, 0010$    & $1111, 0010, 0100$ \\ \hline
$s_{1121}s_{1342}$
& $1111, 1121, 0010$    & $0001, 0010, 1100$ \\ \hline
\end{tabular}}
\label{F4ShortTable}
\end{table}

{\bf The case of {\mathversion{bold}$E_6$}.}
Due to Theorem  \ref{RANKTWO},
given a word $w$
containing a Min-NGT, there exist three roots,
$\alpha_1$, $\alpha_2$, $\alpha_3$, $\alpha_4 $,
with the following properties:

\begin{enumerate}\itemsep=0pt
    \item The roots $\alpha_1$, $\alpha_2$, $\alpha_3$, $\alpha_4$
satisfy the properties of the set of simple roots of
a $D_4$ root system using notation from \cite{Bo}.
    \item The NGT in $\la(w)$ can be written
$\be=\alpha_1 + \alpha_2 + \alpha_3 + \alpha_4$,
$\ga=\alpha_1 + 2\alpha_2 + \alpha_3 + \alpha_4$,
$\al=\alpha_2$.
    \item The nine roots $\alpha_2$, $\alpha_1 + \alpha_2$,
$\alpha_2 + \alpha_3$, $\alpha_2 + \alpha_4$,
$\alpha_1 + \alpha_2 + \alpha_3$, $\alpha_1 + \alpha_2 + \alpha_4$,
$\alpha_2 + \alpha_3 + \alpha_4$,
$\alpha_1 + \alpha_2 + \alpha_3 + \alpha_4$,
$\alpha_1 + 2\alpha_2 + \alpha_3 + \alpha_4$ form the
intersection $\la(w)\cap R^4$ for $R^4\subset R$ generated by
$\{\al_i\}$.
\end{enumerate}

\begin{table}[t]
\caption{$E_6$ Min-NGT: $D_4$-blocks.}
\vspace{1mm}
\centering
{\small
\begin{tabular}{|c|c|c|} \hline
$w$
& NGT
& $\alpha_1$, $\alpha_2$, $\alpha_3$, $\alpha_4$ \\ \hline
$324363243$
& $011101, 012101, 001000$  & $010000, 001000, 000100, 000001$ \\ \hline
$2132436321432$
& $112101, 122101, 010000$  & $100000, 010000, 001000, 001101$ \\ \hline
$4325436324354$
& $012111, 012211, 000100$  & $001000, 000100, 000010, 011001$ \\ \hline
$321432632143263$
& $111101, 112101, 001000$  & $000100, 001000, 000001, 110000$ \\ \hline
$324354632435463$
& $011111, 012111, 001000$  & $010000, 001000, 000001, 000110$ \\ \hline
$3214325436321432543$
& $122211, 123211, 001000$  & $010000, 001000, 000100, 111111$ \\ \hline
$632143254363214325436$
& $123211, 123212, 000001$  & $001000, 000001, 011100, 111110$ \\ \hline
$21324354632143254363243$
& $122211, 123211, 001000$  & $010000, 001000, 000100, 111111$ \\ \hline
$32143254632143254363243$
& $122211, 123211, 001000$  & $010000, 001000, 000100, 111111$ \\ \hline
$32143263214325436321432$
& $112111, 122111, 010000$  & $100000, 010000, 001000, 001111$ \\ \hline
$32435463214325436324354$
& $112111, 112211, 000100$  & $001000, 000100, 000010, 111001$ \\ \hline
$43215432632143254363243$
& $122211, 123211, 001000$  & $010000, 001000, 000100, 111111$ \\ \hline
$2132435463214325436324354$
& $122111, 122211, 000100$  & $000010, 000100, 011000, 111001$ \\ \hline
$4321543263214325436321432$
& $112211, 122211, 010000$  & $100000, 010000, 001100, 001111$ \\ \hline
$32143254363214325436321432543$
& $111111, 112111, 001000$  & $000001, 001000, 110000, 000110$ \\ \hline
\end{tabular}}
\label{E6Table}
\end{table}

The root $\alpha_2$ is always a simple root of $E_6$
but $\alpha_1$, $\alpha_3$ and $\alpha_4$ need not be.
The three roots~$\alpha_1$,~$\alpha_3$ and
$\alpha_4$ will not be elements of $\lambda(w)$.
Together with the nine roots above these are all twelve
positive roots in the $D_4$ subsystem $R^4$ determined by
$\alpha_1$, $\alpha_2$, $\alpha_3$, $\alpha_4$,
i.e., the condition {\em nine but not twelve}
from (\ref{rankthreed}) is satisf\/ied.

Table~\ref{E6Table} shows each of the Min-NGT words in $E_6$.
Also included are the roots that make up the NGT and a
choice of $\alpha_1$, $\alpha_2$, $\alpha_3$, $\alpha_4$ def\/ining
an appropriate $D_4$ subsystem. Since all of the words
can not be written as compositions of pairwise orthogonal
ref\/lections, we uniformly put them in the {\em shortlex form}:
the lexicographical ordering from left to right, but with
the composition from right to left.
We use a one-line representation
of the roots from $E_6$ where the coef\/f\/icient
of the exceptional simple root is placed {\em the last},
i.e.,
$\begin{array}{c}
abcde \\
f
\end{array}$ is written $abcdef$.

\subsection*{Acknowledgements}
Partially supported by NSF grant DMS--0800642.

\pdfbookmark[1]{References}{ref}
\LastPageEnding

\end{document}